\documentclass[reqno]{amsart}
\usepackage{amssymb,amscd}

\begin{document}

\let\kappa=\varkappa
\let\eps=\varepsilon
\let\phi=\varphi

\def\Z{\mathbb Z}
\def\R{\mathbb R}
\def\C{\mathbb C}
\def\Q{\mathbb Q}

\def\OO{\mathcal O}
\def\CP{\C{\mathrm P}}
\def\RP{\R{\mathrm P}}
\def\conj{\overline}
\def\Beta{\mathrm{B}}
\def\p{\partial}

\renewcommand{\Im}{\mathop{\mathrm{Im}}\nolimits}
\renewcommand{\Re}{\mathop{\mathrm{Re}}\nolimits}
\newcommand{\codim}{\mathop{\mathrm{codim}}\nolimits}
\newcommand{\id}{\mathop{\mathrm{id}}\nolimits}
\newcommand{\Aut}{\mathop{\mathrm{Aut}}\nolimits}

\newtheorem*{thm}{Theorem}
\newtheorem*{cor}{Corollary}

\theoremstyle{definition}
\newtheorem*{rem}{Remark}

\title{Note on the paper of Fu and Wong on strictly pseudoconvex domains with K\"ahler--Einstein
Bergman metrics}
\author{Stefan Nemirovski}
\address{%
Steklov Mathematical Institute, 119991 Moscow, Russia\hfill\break
\phantom{hh} \& \hfill\break
\phantom{hf} Ruhr-Universit\"at Bochum, D-44780 Bochum, Germany}
\email{stefan@mi.ras.ru}
\author{Rasul Shafikov}
\address{Department of Mathematics, the University of Western Ontario,
London, Canada  N6A 5B7}
\email{shafikov@uwo.ca}
\begin{abstract}
It is shown that the Ramadanov conjecture implies the Cheng conjecture. 
In particular it follows that the Cheng conjecture holds in dimension two.
\end{abstract}

\maketitle

In this brief note we use our uniformization result from~\cite{NS1,NS2}
to extend the work of Fu and Wong~\cite{FW} on the relationship
between two long-standing conjectures about the behaviour of the
Bergman metric of a strictly pseudoconvex domain in $\C^n$, $n\ge 2$.

\smallskip
\noindent
{\bf 1.}
Let $D\Subset\C^n$ be an arbitrary bounded domain.
The {\it Bergman kernel function\/} of~$D$ can be defined by the formula
$$
K_D(z):=\sum_{j=1}^\infty \phi_j(z)\conj{\phi_j(z)},
$$
where $\{\phi_j\}_{j=1,\ldots,\infty}$ is any orthonormal basis
of the space $L^2\OO(D)$ of square-integrable holomorphic functions
in~$D$.

It is a standard result that the function $\log K_D(z)$ is strictly plurisubharmonic
and the positive $(1,1)$-form
$$
k_D:= i\p\conj\p \log K_D(z)
$$
is invariant with respect to biholomorphic mappings between bounded domains.
The {\it Bergman metric\/} on~$D$ is the K\"ahler metric associated
with this K\"ahler form.

\smallskip
\noindent
{\bf 2.}
Fefferman~\cite{Fe} (see also~\cite{BdMS}) established the following
deep result on the boundary behaviour of the kernel function of a
smoothly bounded strictly pseudoconvex domain~$D$. Let
$\rho\in C^\infty (\conj D)$ be a defining function for~$D$. Then
there is a decomposition
\begin{equation}
\label{Feff}
K_D(z)=\phi(z) \rho(z)^{-(n+1)} + \psi(z) \log |\rho(z)|
\end{equation}
where the functions $\phi,\psi\in C^\infty(\conj D)$, and $\phi\ne 0$
everywhere on~$\p D$. Note that the latter property implies that
the Bergman metric of a strictly pseudoconvex domain is complete.

Although the kernel function is defined globally, its asymptotic behaviour
as $z\to z_0\in \p D$ depends only on the local CR geometry of the boundary
at the point~$z_0$ (see~\cite{Fe}). For instance, the kernel function of
the unit ball $B\subset\C^n$ can be explicitly decomposed as
$$
K_B(z)=\frac{n!}{\pi^n}\,(1-\|z\|^2)^{-(n+1)}
$$
with identically vanishing logarithmic term. Thus the coefficient~$\psi$
in the logarithmic term of~(\ref{Feff}) vanishes to infinite order at the boundary
for any strictly pseudoconvex domain whose boundary is spherical (i.\,e.,
locally CR diffeomorphic to the unit sphere~$\p B\subset\C^n$).

The {\it Ramadanov conjecture\/}~\cite{Ra} asserts that, conversely,
the vanishing condition $\psi=O(\rho^\infty)$ implies that the
boundary of~$D$ is spherical. This conjecture has been proved
for domains in~$\C^2$ by Graham and Burns~\cite{Grm} and Boutet de Monvel~\cite{BdM}.

\smallskip
\noindent
{\bf 3.}
A classical problem, proposed in different forms by Bergman, Hua, and Yau,
asks to describe the domain in terms of the differential-geometric
properties of its Bergman metric. For example, a well-known theorem
of Lu Qi-Keng~\cite{Lu} states that a bounded domain with complete Bergman
metric of constant holomorphic sectional curvature is biholomorphic
to the ball.

The {\it Cheng conjecture\/}~\cite{Che}\footnote{%
The formulation in~\cite{Che} is somewhat vague; the precise
statement below is taken from~\cite{FW}.} asserts that the
hypotheses of Lu's theorem can be weakened for a smoothly bounded
strictly pseudoconvex domain. Namely, such a domain has to be
biholomorphic to the ball if and only if its Bergman metric
is K\"ahler--Einstein.

\begin{thm}
\label{ChRa}
The Cheng conjecture in\/~$\C^n$ follows from the Ramadanov conjecture in\/~$\C^n$.
\end{thm}

Since the Ramadanov conjecture is known to be true in dimension~$2$,
we obtain the following result.

\begin{cor}
The Bergman metric of a smoothly bounded strictly pseudoconvex domain $D\Subset\C^2$
is K\"ahler--Einstein if and only if this domain is biholomorphic to the ball.
\end{cor}

\begin{rem}
Fu and Wong~\cite{FW} proved these results for {\it simply connected\/}
domains using a weaker uniformization result of Chern and Ji~\cite{cj}
and stated the general case as an open question.
\end{rem}

\begin{proof}[Proof of the theorem]
Suppose that the Bergman metric of a strictly pseudoconvex domain~$D$
is K\"ahler--Einstein. Fu and Wong~\cite{FW} computed (rather ingeniously)
that in this case the logarithmic coefficient $\psi$ in the decomposition~(\ref{Feff})
vanishes to infinite order at~$\p D$. Assuming the Ramadanov conjecture,
we conclude that the boundary of $D$ is spherical. Hence, by the
uniformization theorem (see \cite[Thm.~A.2]{NS1} and~\cite[Cor.~3.2]{NS2}) the domain $D$
is covered by the unit ball.

Since $D$ is the quotient of the ball by the group of holomorphic
deck transformations, there is a natural complete metric of
constant holomorphic sectional curvature on $D$ obtained by 
taking the quotient of the standard invariant metric on the ball.
According  to Cheng and Yau~\cite{CY}, the complete K\"ahler--Einstein metric
on $D$ is unique up to a constant factor.  Hence, the Bergman metric of~$D$ 
is proportional to the quotient metric and so has constant holomorphic sectional 
curvature. It follows that $D$ is biholomorphic to the ball by Lu's theorem~\cite{Lu}.
\end{proof}

\end{document}